\documentclass[12pt]{article}
\usepackage{amsmath, amsthm, amssymb,graphicx}
\numberwithin{equation}{section}
\begin{document}
\author{Ajai Choudhry}
\title{Cyclic pentagons and hexagons with integer sides, diagonals and areas}
\date{}
\maketitle

 \begin{abstract}
In this paper we obtain  cyclic pentagons and hexagons with rational sides, diagonals and area all of which are expressed  in terms of  rational functions of several arbitrary  rational parameters. On suitable scaling, we obtain cyclic pentagons and hexagons whose sides, diagonals and area are all given by integers.  

\end{abstract}



\section{Introduction}
A polygon is called a cyclic polygon if all its vertices lie on the circumference of a circle. A polygon all of whose  sides, diagonals and the area are integers is called a Heron polygon. A cyclic Heron polygon  is called a Brahmagupta polygon. Further, a polygon with rational sides, diagonals  and rational area is called a rational polygon. This paper is concerned with Brahmagupta  pentagons and hexagons, that is, 
cyclic pentagons and hexagons with integer sides, diagonals and areas.

It would be recalled that the area $K$ of a triangle is given by Heron's  formula,
\begin{equation}
K=\sqrt{s(s-a)(s-b)(s-c)} \label{heron}
 \end{equation}
where $a,\,b,\,c$ are the sides of the triangle and $s$ is the semi-perimeter given by $s=(a+b+c)/2$.
Further, if  $a_1,\,a_2,\,a_3,\,a_4$ are the consecutive sides of a cyclic quadrilateral, the formulae for the area $K$ and the two diagonals $d_1$ and $d_2$, given by Brahmagupta \cite[p. 187]{Ev} in the seventh century A. D., are as follows:
\begin{equation}
\begin{aligned}
K&=\sqrt{(s-a_1)(s-a_2)(s-a_3)(s-a_4)},\\
d_1&=\sqrt{(a_1a_2+a_3a_4)(a_1a_3+a_2a_4)/(a_1a_4+a_2a_3)},\\
d_2&=\sqrt{(a_1a_3+a_2a_4)(a_1a_4+a_2a_3)/(a_1a_2+a_3a_4)}, 
\end{aligned}
\label{Brahmaformulae}
\end{equation}
where $s$ is the semi-perimeter, that is,
\begin{equation}
s=(a_1+a_2+a_3+a_4)/2. \label{vals}
\end{equation}
Euler gave a  complete parametrisation of all rational triangles (see Dickson \cite[p. 193]{Di}). 
More recently, Sastry \cite{Sa1} has given  a complete parametrisation of all cyclic quadrilaterals whose sides, diagonals, circumradius and area are all given by integers.

The formulae for the area of cyclic pentagons and cyclic hexagons are very complicated and were discovered only towards the end of the twentieth century   by Robbins (\cite{Ro1}, \cite{Ro2}). In fact, if $K$ is the area of a cyclic pentagon/ hexagon, then $K^2$ satisfies an equation of degree 7 whose coefficients are symmetric functions of the sides of the cyclic pentagon/ hexagon. Similarly, the formulae 
for the diagonals of cyclic pentagons and hexagons are given by equations of degree 7 \cite[p. 37 and p. 41]{BM}.
In view of these complicated formulae, the areas and the diagonals
of cyclic pentagons and hexagons with rational sides are, in general, not rational.  There has thus been considerable interest in constructing cyclic pentagons and hexagons whose sides, diagonals and the area are all given by rational numbers. 

 Euler \cite[p. 221]{Di}) and Sastry \cite{Sa2} have independently described  methods of constructing   cyclic polygons with five or more sides such that all the sides, the diagonals and the area of the polygon are rational. Recently, Bucchholz and MacDougall \cite{BM} have carried out computational searches for cyclic  rational  pentagons and hexagons.  While these computational searches yielded several cyclic rational pentagons, Bucchholz and MacDougall did not find even a single hexagon with rational sides and area and with all nine diagonals rational. Till date a  parametrisation of cyclic pentagons  and hexagons does not appear to have been published.

We note that geometric figures such as triangles, quadrilaterals and polygons with rational sides, diagonals (where applicable), and areas yield, on appropriate scaling, similar geometric figures with integer sides, diagonals (where applicable), and areas. It therefore suffices to obtain cyclic pentagons and hexagons with rational sides, diagonals and areas.

In this paper we obtain a cyclic pentagon with rational sides, diagonals and area all of which are expressed  in terms of  rational functions of several arbitrary  rational parameters. We also show how additional cyclic rational pentagons may be obtained. The complete parametrisation of all  cyclic rational pentagons is given by a finite number of such parametrisations.   We also give a parametrisation of cyclic hexagons whose sides, the nine diagonals and the area are all rational.

\section{Parametrisations of rational  triangles and  quadrilaterals}

We will construct a cyclic pentagon by juxtaposing a rational  triangle and a cyclic rational  quadrilateral having the same circumradius. Similarly we will construct a cyclic  hexagon by juxtaposing two cyclic rational  quadrilaterals having the same circumradius. Accordingly, we give below the already known complete parametrisation of rational   triangles and  quadrilaterals.

The complete parametrisation of  an arbitrary rational  triangle was given by Euler (see Dickson \cite[p. 193]{Di}) who proved that the sides $a,\,b,\,c$ of such a triangle   are given by
\begin{equation}
\begin{aligned}
a &= k(m^2+n^2)/(mn),\\
b &= k(p^2+q^2)/(pq),\\
 c&=k(pn+qm)(pm-qn)/(pqmn).
\end{aligned}
\label{euler}
\end{equation}
where $m,\,n,\,p,\,q$ and $k$ are arbitrary  rational parameters. The area of the triangle is $k^2(pn+qm)(pm-qn)/(pqmn)$.

Further,  the circumradius $R$ of a   triangle whose sides are  $a,\,b,\,c$ is given by
\begin{equation}
R=abc/\{4\sqrt{s(s-a)(s-b)(s-c)}\},
\label{Rtrg}
\end{equation}
where $s$ is the semi-perimeter, that is, $s=(a+b+c)/2$. Thus, the circumradius $R$ is given in terms of the sides $a,\,b,\,c$ explicitly by the formula,
\begin{equation}
R=abc/\sqrt{(a+b+c)(a+b-c)(b+c-a)(c+a-b)},
\label{Rtrgabc}
\end{equation}
 It follows that the circumradius of a triangle with sides \eqref{euler} is 
\begin{equation}
 k(m^2+n^2)(p^2+q^2)/(4mnpq). \label{circumradiuseuler}
\end{equation}

The complete parametrisation of  an arbitrary cyclic quadrilateral with rational sides, diagonals and area was obtained by Sastry \cite[p. 170]{Sa1}  who proved that the sides $a_1,\,a_2,\,a_3,\,a_4$, the diagonals $d_1,\,d_2$, the area $K$ and the circumradius $R$  of such a quadrilateral   are given by
\begin{equation}
\begin{aligned}
a_1&=\{t(u + v) + ( 1 - uv )\}\{u + v - t( 1 - uv )\},\\
a_2&=( 1 + u^2 ) ( v - t )( 1 + t v ),\\
a_3 &= t( 1 + u^2 )( 1 + v^2 ), \\
a_4&=( 1 + v^2 )( u - t )( 1 + tu ),\\
d_1&=u( 1 + t^2 )( 1 + v^2),\\
d_2&=v( 1 + t^2 )( 1 + u^2),\\
K&=uv\{2t(1- uv) - (u + v)(1 - t^2 )\}\\
& \quad \quad \quad \times \{2(u + v)t + (1 - uv)(1 - t^2)\},\\
R&=( 1 + t^2 )( 1 + u^2)(1+v^2)/4.
\end{aligned}
\label{sastry}
\end{equation}
where $t,\,u$ and $v$  are arbitrary  rational parameters. 

We also note that  the circumradius $R$ of a cyclic quadrilateral with sides $a_1,\,a_2,\,a_3$ and $a_4$ is given by the following formula of Paramesvara \cite{Gu1}:
\begin{equation}
R=\frac{\displaystyle 1}{\displaystyle 4}\ \sqrt{\frac{\displaystyle(a_1a_2+a_3a_4)(a_1a_3+a_2a_4)/(a_1a_4+a_2a_3)}{\displaystyle (s-a_1)(s-a_2)(s-a_3)(s-a_4)}},
\end{equation}
where $s$ is the semi-perimeter of the quadrilateral, that is, $s=(a_1+a_2+a_3+a_4)/2$. Thus, the circumradius $R$ is given in terms of the sides $a_1,\,a_2,\,a_3,\,a_4$ explicitly by the formula,
\begin{multline}
R=\{(a_1a_2+a_3a_4)(a_1a_3+a_2a_4)(a_1a_4+a_2a_3)\}^{1/2}\{(-a_1+a_2+a_3+a_4)\\
\times (a_1-a_2+a_3+a_4)(a_1+a_2-a_3+a_4)(a_1+a_2+a_3-a_4)\} ^{-1/2}.
\label{Param}
\end{multline}

\section{Cyclic pentagons with integer sides, \\diagonals and area}
Any diagonal of a  pentagon divides it into two parts: a triangle and a quadrilateral with a common side which is one of the diagonals of the pentagon. Accordingly we may  consider a cyclic pentagon $ABCDE$ (see Figure 1) with rational sides and diagonals as  being made up by juxtaposing the cyclic rational quadrilateral  $ABCE$ and the rational triangle $CDE$. The lengths of the sides of the quadrilateral $ABCE$ and the triangle $CDE$ can be chosen in several ways. We consider these possibilities in the next two subsections.
\begin{figure}[h]
\includegraphics[scale =1.5, trim=0 3cm 0 6cm, clip, width=\textwidth]{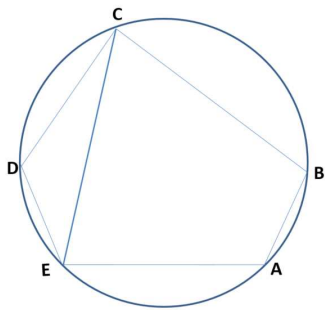}
\caption{Cyclic pentagon}
\centering
\end{figure}

\subsection{}
In this subsection, we take the sides $AB, \; BC,\;CE,\;$ and $EA$ as having the lengths $a_1,\,a_2,\,a_3,\,a_4$ respectively as given by \eqref{sastry} while the lengths of the three sides of the triangle $CDE$ may be taken as given by \eqref{euler}. 

Since the  quadrilateral  $ABCE$ and the triangle $CDE$ are inscribed in the same circle, the circumradii of the  quadrilateral  $ABCE$ and the  triangle $CDE$ are equal, and hence, on using the formulae \eqref{circumradiuseuler} and \eqref{sastry}, we get,
\begin{equation}
k(m^2+n^2)(p^2+q^2)/(4mnpq) = (t^2+1)(v^2+1)(u^2+1)/4. \label{condR}
\end{equation}

Now $CE$ is also a side of the triangle $CDE$ and therefore, its length may be taken, without loss of generality, either as $a$ or  $c$ as given by \eqref{euler}.

We first consider the case when $CE$ has length $a$ as given by \eqref{euler}. Equating the two lengths of $CE$, we get  
\begin{equation}
k(m^2+n^2)/(mn) = t(u^2+1)(v^2+1). \label{condCE}
\end{equation}

Eqs.~\eqref{condR} and \eqref{condCE} are readily solved and we get the following two solutions:
\begin{equation}
t=p/q,\quad k=pmn(u^2+1)(v^2+1)/\{q(m^2+n^2)\}, \label{sol1}
\end{equation}
and
\begin{equation}
t=q/p,\quad k=qmn(u^2+1)(v^2+1)/\{p(m^2+n^2)\}. \label{sol2}
\end{equation}

Using the solution \eqref{sol1} and the formulae \eqref{euler} and \eqref{sastry}, we readily find the sides of the pentagon. Denoting the lengths of the sides $AB, \,BC,\,CD,$  $DE$ and $EA$ by $s_1,\,s_2,\ldots,\,s_5$ respectively, we get,
\begin{equation}
\begin{aligned}
s_1 &= (-quv+pu+pv+q)(puv+qu+qv-p)/q^2, \\
s_2 &= (u^2+1)(qv-p)(pv+q)/q^2, \\
s_3 &= mn(u^2+1)(v^2+1)(p^2+q^2)/\{q^2(m^2+n^2)\},\\
 s_4 &= (u^2+1)(v^2+1)(mq+np)(mp-nq)/\{q^2(m^2+n^2)\},\\
 s_5 &= (v^2+1)(qu-p)(pu+q)/q^2,
\end{aligned}
\label{sidespent1}
\end{equation}
where $m,\,n,\,p,\,q,\,u$ and $v$ are arbitrary rational parameters.

The circumradius $R$ of the pentagon, given by either side of \eqref{condR}, may be written explicitly as follows:
\begin{equation}
R= (u^2+1)(v^2+1)(p^2+q^2)/(4q^2). \label{Rpent1}
\end{equation}

The area of the pentagon is rational and is  easily worked out, being the sum of the rational areas of the triangle $CDE$ and the quadrilateral $ABCE$.

As regards the diagonals, the length of the diagonal $CE$, which is the common  side of the triangle $CDE$ and the quadrilateral $ABCE$
diagonal $CE$, is readily found and is given by
\[
CE= p(u^2+1)(v^2+1)/q.\]
The two diagonals $AC$ and $BE$, being diagonals of the cyclic quadrilateral $ABCE$, may be  found using the formulae \eqref{sastry} and are given by 
\begin{equation}
\begin{aligned}
AC = u(p^2+q^2)(v^2+1)/q^2, \quad BE = v(p^2+q^2)(u^2+1)/q^2.
\end{aligned}
\end{equation}
To find the length $d$ of the diagonal  $AD$, we note that the circumradius  of the triangle $ADE$ is given by \eqref{Rpent1}, and since the sides of the triangle $ADE$ have lengths $d,\,s_4$ and $s_5$, in view of the formula \eqref{Rtrgabc}, we must have 
\begin{multline}
ds_4s_5/\sqrt{(d+s_4+s_5)(d+s_4-s_5)(s_4+s_5-d)(s_5+d-s_4)}\\
=(u^2+1)(v^2+1)(p^2+q^2)/(4q^2), \label{condd1}
\end{multline}
where $s_4,\,s_5$ are given by \eqref{sidespent1}. 

We also note that the circumradius  of the cyclic quadrilateral  $ABCD$ is given by \eqref{Rpent1}. Since the sides of the quadrilateral $ABCD$ have lengths $s_1,\,s_2,\,s_3$ and $d$, it follows from the formula \eqref{Param} that
\begin{multline}
\{(s_1s_2+s_3d)(s_1s_3+s_2d)(s_1d+s_2s_3)\}^{1/2}\{(-s_1+s_2+s_3+d)\\
\times (s_1-s_2+s_3+d)(s_1+s_2-s_3+d)(s_1+s_2+s_3-d)\} ^{-1/2}\\
=(u^2+1)(v^2+1)(p^2+q^2)/(4q^2), \label{condd2}
\end{multline}
where $s_1,\,s_2,\,s_3$ are given by \eqref{sidespent1}. 

Each of the two equations \eqref{condd1} and \eqref{condd2} is satisfied by four rational values of $d$ but there is  one common root which gives the length $d$ of the diagonal $AD$. We thus get
\[
AD=(p^2+q^2)(v^2+1)(mu-n)(nu+m)/\{(m^2+n^2)q^2\}.
\]

Finally, the  diagonal $BD$ may be obtained just as the diagonal $AD$ by solving two equations obtained by equating the circumradius of the pentagon first to the circumradius of the triangle $BCD$  and then to the circumradius  of the quadrilateral  $ABDE$ and finding the common root. We thus get
\[
BD=(u^2+1)(mqv+npv-mp+nq)(mpv-nqv+mq+np)/\{(m^2+n^2)q^2\}.
\]

Thus, the cyclic pentagon $ABCDE$ with rational sides $s_i,\;i=1,\,\ldots,\,5$, given by \eqref{sidespent1} in terms of arbitrary  rational parameters $m,\,n,\,p, \,q,\,u$ and $v$  has rational circumradius \eqref{Rpent1}, rational area and rational diagonals as seen above.

We note that we have obtained the above cyclic pentagon just by equating the circumradii and the side $CE$ of the triangle $CDE$ and the quadrilateral $ABDE$. We also need to ensure that the triangle and the quadrilateral can actually  be constructed  on opposite sides of $CE$. 

 If for certain values of the parameters, the common side $CE$ is twice the circumradius, then $CE$ is a diameter of the circumcircle and if $s_i,\;i=1,\,\ldots,\,5$, are all positive,  the triangle $CDE$ and the quadrilateral $ABDE$ can always be constructed on opposite sides of the diameter $CE$.  These conditions are satisfied when $p=q,\;u > 1,\; m > n$ and $(u+1)/(u-1) > v > 1$. As a numerical example, when $p=q=1,\;u=3,\; v=3/2,\; m=2,\; n=1$, we get a cyclic pentagon with sides $AB, \, BC, \, CD, \, DE,\,EA$ being $8,\, 25/2,\,  26,\,  39/2,\, 26$, the diagonals have lengths $39/2,\;  65/2,\;63/2,\; 30,\;  65/2$, the circumradius has length 65/4 and the area of the pentagon is 537.

When the side $CE$ divides the circle into two unequal segments,  there may exist  certain values of the parameters for which  it is possible that both the triangle and the quadrilateral can  be constructed only in the same  segment and, in such a case, the desired  cyclic pentagon cannot be constructed. To ensure that we actually get a cyclic pentagon, we observe that there is no loss of generality in   taking the triangle to be in the minor segment  and the quadrilateral in the major segment. 

If we choose the parameters $m,\,n,\,p,\,q$ satisfying the conditions 
\[
q/p > {\rm min}(u,v) > p/q > (uv-1)/(u+v),\quad {\rm and \quad } qm > pn,
\]
then it is easily seen that the sides $CD$ and $DE$ are both shorter than $CE$ and hence  the triangle is necessarily in the minor segment of the circumcircle. Further, one of the diagonals of the quadrilateral is longer than $CE$  and so the quadrilateral can be constructed in the major segment, and we will actually have a cyclic pentagon with the desired properties.

 As an example, taking $m=3,\;n=1,\;p=2,\;q=3,\;u=1,\; v=2$, we get a cyclic pentagon whose sides are $AB, \, BC, \, CD, \, DE,\,EA$  are $11/3,\; 56/9,\; 13/3,\; 11/3,\; 25/9$ respectively, the diagonal $CE$ has length 20/3, the remaining diagonals have lengths $65/9,\; 52/9,\; 323/45, \; 52/9$, the circumradius is 65/18 and the area of the pentagon is 28.

 We note that using the solution \eqref{sol2} of Eqs.~\eqref{condR} and \eqref{condCE}, we also get a cyclic pentagon but this is essentially equivalent to the pentagon that we have already obtained using the first solution \eqref{sol1}.

We now consider the case when the side $CE$ of the triangle $CDE$ is taken as  $c$ given by \eqref{euler} while the sides $AB, \; BC,\;CE,\;$ and $EA$ of the quadrilateral are, as before, having lengths $a_1,\,a_2,\,a_3,\,a_4$ respectively  given by \eqref{sastry}. Proceeding as in subsection 3.1, we find the sides $s_i,\;i=1,\,\ldots,\,5$ and the circumradius of the pentagon. It turns out, however, that by a suitable transformation of the parameters, this solution is seen to be equivalent to the solution already obtained above.

\subsection{} In the  previous subsection, we have considered $CE$ as a side of the quadrilateral $ACDE$ as well as a side of the triangle $CDE$ and equated the two lengths.  While we had taken the length of the side $CE$ of the quadrilateral $ABCE$ as $a_3$, we can also choose the side $CE$ to have  length   $a_1$ or $a_2$ or $a_4$ as given by \eqref{sastry}. Since the lengths $a_2$ and $a_4$ get interchanged if we interchange the parameters $u$ and $v$,  there is no loss of generality in taking the length of $CE$ either as $a_1$ or $a_2$. If we consider $CE$ as one of the sides of the triangle $CDE$, we have already noted that without loss of generality,  we may take its length either as $a$ or $c$ as given by \eqref{euler}. 

We thus get four cases according to the choice of lengths of the side $CE$. If we now proceed as in subsection 3.1 and impose the two conditions that the circumradii of the triangle $CDE$ and the quadrilateral $ABCE$ are equal, and the two lengths of   $CE$ are also equal, it is noteworthy that in each case, the two conditions are readily solvable, and further, in each case, the five diagonals of the pentagon are also rational, and their lengths can be computed as in subsection 3.1. 

We thus get a finite number of   parametrisations  for a cyclic pentagon with rational sides, diagonals and area. All such cyclic pentagons arise from one of these parametrisations. We do not give these remaining parametrisations explicitly as any interested reader can readily find them by following the method described in subsection 3.1. 

\section{Cyclic hexagons with integer sides, \\diagonals and area} 
We will now construct cyclic hexagons by juxtaposing two cyclic rational quadrilaterals which  have the same circumradius and which have a common side which becomes a diagonal of the hexagon (see Figure 2).
\begin{figure}[h]
\includegraphics[scale=1.5, trim=0 6cm 0 3cm, clip,  width=\textwidth]{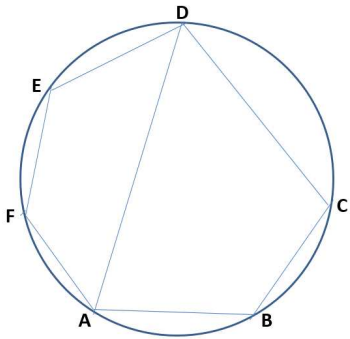}
\caption{Cyclic hexagon}
\centering
\end{figure}

It has been shown in \cite{BM} that  merely juxtaposing  two cyclic rational quadrilaterals with the same circumradius and a common side 
does not generally  yield a cyclic hexagon with rational diagonals. We will, however, choose the sides of the two quadrilaterals in such a manner that we are able to obtain a cyclic hexagon whose sides, diagonals, circumradius and area are all rational.

 We will take the sides of the cyclic quadrilateral $ABCD$ as being given by Sastry's formulae \eqref{sastry} with the parameters $t,\,u,\,v,$ having values $t,\,u_1$ and $v_1$ respectively. We thus get the length of the circumradius of the quadrilateral $ABCD$ as
$(t^2+1)(u_1^2+1)(v_1^2+1)/4$.  We now  specifically take the side $AD$ as having  length $a_3$, that is, we take,
$AD=t(u_1^2+1)(v_1^2+1)$.

For the second cyclic quadrilateral $ADEF$, we take the sides as being given by the  formulae \eqref{sastry} with the parameters $t,\,u,\,v,$ having values $t,\,u_2$ and $v_2$ respectively. We thus get the length of the circumradius of the quadrilateral $AFEDD$ as
$(t^2+1)(u_2^2+1)(v_2^2+1)/4$.  Further, we again take the side $AD$ as having  length $a_3$ given by \eqref{sastry}, so that we get
$AD=t(u_2^2+1)(v_2^2+1)$.

It is now readily seen that the two conditions, that both the cyclic quadilaterals have the same circumradius, and that they have a side in common, are simultaneously satisfied if the parameters $u_i,\,v_i$, satisfy the following condition:
\begin{equation}
(u_1^2+1)(v_1^2+1)=(u_2^2+1)(v_2^2+1). \label{condhex}
\end{equation}

Now Eq.~\eqref{condhex} may be written as,
\begin{equation}
\begin{aligned}
(u_1v_1+1)^2+(u_1-v_1)^2 &= (u_2v_2+1)^2+(u_2-v_2)^2,\\
{\rm or,} \quad \quad (u_1v_1+1)^2-(u_2v_2+1)^2&=(u_2-v_2)^2-(u_1-v_1)^2,\\
{\rm or,} \quad \quad (u_1v_1+u_2v_2+2)(u_1v_1-u_2v_2)& = -(u_1+u_2-v_1-v_2)\\
& \quad \quad \times (u_1-u_2-v_1+v_2),
\end{aligned}
\end{equation}
and is therefore equivalent to the following two linear equations in $v_1,\,v_2$:
\begin{equation}
\begin{aligned}
u_1v_1+u_2v_2+2 = -m(u_1+u_2-v_1-v_2),\\
m(u_1v_1-u_2v_2) = u_1-u_2-v_1+v_2.
\end{aligned}
\end{equation}
where $m$ is an arbitrary rational parameter. We thus get the following solution of Eq.~\eqref{condhex}:
\begin{equation}
\begin{aligned}
v_1&= \{(u_1+u_2)u_2m^2+(2u_1+2u_2)m-u_1u_2+u_2^2+2\}\\
& \quad \quad \times \{(u_1+u_2)m^2+(-2u_1u_2+2)m-u_1-u_2\}^{-1},\\
 v_2& = \{(u_1+u_2)u_1m^2+(2u_1+2u_2)m+u_1^2-u_1u_2+2\}\\
& \quad \quad \times \{(u_1+u_2)m^2+(-2u_1u_2+2)m-u_1-u_2\}^{-1}. \label{hexv12}
\end{aligned}
\end{equation}
 
With these values of $v_1,\,v_2$, we can readily work out the sides of the two cyclic quadilaterals $ABCD$ and $ADEF$ as well as their common circumradius in terms of arbitrary rational parameters $m,\,t,\,u_1$ and $u_2$. We will denote  the lengths of the sides $AB,\,BC,\,CD,\,DE$ and $EF$ of the hexagon $ABCDEF$ by $h_i(m,\,t,\,u_1,\,u_2),\;i=1,\ldots,\,6$, respectively. On appropriate scaling, we may write the lengths of the sides  of the hexagon $ABCDEF$ as follows:
\begin{equation}
\begin{aligned}
h_1(m,\,t,\,u_1,\,u_2)&=\{(u_1+u_2)((u_1+u_2)t-u_1u_2+1)m^2\\
& \quad \quad +((-2u_1^2u_2+4u_1+2u_2)t-2u_1^2-4u_1u_2+2)m\\
& \quad \quad +(-u_1^2-2u_1u_2+u_2^2+2)t+u_1^2u_2-u_1u_2^2-3u_1-u_2\}\\
& \quad \quad \times \{(u_1+u_2)((u_1u_2-1)t+u_1+u_2)m^2\\
& \quad \quad +((2u_1^2+4u_1u_2-2)t-2u_1^2u_2+4u_1+2u_2)m\\
& \quad \quad +(-u_1^2u_2+u_1u_2^2+3u_1+u_2)t-u_1^2-2u_1u_2+u_2^2+2\},\\
h_2(m,\,t,\,u_1,\,u_2)&=-(u_1^2+1)\{(u_1+u_2)(t-u_2)m^2+((-2u_1u_2+2)t\\
& \quad \quad -2u_1-2u_2)m+(-u_1-u_2)t+u_1u_2-u_2^2-2\}\\
& \quad \quad \times \{(u_1+u_2)(u_2t+1)m^2 +((2u_1+2u_2)t\\
& \quad \quad -2u_1u_2+2)m+(-u_1u_2+u_2^2+2)t-u_1-u_2\},\\
h_3(m,\,t,\,u_1,\,u_2)&=(u_2^2+1)(m^2+1)(u_1-t)(tu_1+1) \{(u_1+u_2)^2m^2\\
& \quad \quad +4(u_1+u_2)m+u_1^2-2u_1u_2+u_2^2+4\},\\
h_4(m,\,t,\,u_1,\,u_2)&=h_1(m,\,t,\,u_2,\,u_1),\\
h_5(m,\,t,\,u_1,\,u_2)&=h_2(m,\,t,\,u_2,\,u_1),\\
h_6(m,\,t,\,u_1,\,u_2)&=h_3(m,\,t,\,u_2,\,u_1).
\end{aligned}
\label{hex}
\end{equation}

The circumradius $R$ of the hexagon is given by
\begin{multline}
R=(u_1^2+1)(u_2^2+1)(t^2+1)(m^2+1)\\
\times \{(u_1+u_2)^2m^2+(4u_1+4u_2)m+u_1^2-2u_1u_2+u_2^2+4\}/4. \label{Rhex}
\end{multline}

The area of the hexagon, being the sum of the  rational areas of the two cyclic quadrilaterals $ABCD$ and $ADEF$ is naturally rational, and is easily determined.

We will now determine all the 9 diagonals of the hexagon $ABCDEF$. We first find the  three central  diagonals $AD,\,BE,\,CF$, each of which divides the hexagon into two cyclic quadrilaterals. We have, in fact, already found the diagonal $AD$ since it is the common side of the two quadrilaterals $ABCD$ and $ADEF$, and its length is given by
\begin{multline}
AD=t(u_1^2+1)(u_2^2+1)(m^2+1)\\
\times  \{(u_1+u_2)^2m^2+4(u_1+u_2)m+u_1^2-2u_1u_2+u_2^2+4\}.
\end{multline}

To find the length $d$ of the central diagonal $BE$, we note that $d$ is one of the sides of each of the two quadrilaterals  $BCDE$ and $ABEF$ both of which have circumradius $R$. Thus $d$ must satisfy the condition that the circumradius of the quadrilateral $BCDE$ is $R$. Writing the sides $BC,\,CD,\,DE$ simply as $h_2, \; h_3,\,h_4$, and using the formula \eqref{Param}, we get the condition,
\begin{multline}
R=\{(dh_2+h_3h_4)(dh_3+h_2h_4)(dh_4+h_2h_3)\}^{1/2}\{(-d+h_2+h_3+h_4)\\
\times (d-h_2+h_3+h_4)(d+h_2-h_3+h_4)(d+h_2+h_3-h_4)\} ^{-1/2}.
\label{Paramhex1}
\end{multline}
The values of $h_2, \; h_3,\,h_4$ and $R$ are given by \eqref{hex} and \eqref{Rhex}, and on solving Eq.~\eqref{Paramhex1} for $d$, we get four rational roots. Similarly, $d$  must satisfy the condition that the circumradius of the quadrilateral $ABEF$ is $R$, and this condition also results in an equation that is satisfied by four rational values of $d$. The common root gives the length $d$ and we thus get,
\begin{multline}
BE=\{(u_1+u_2)^2m^2+4(u_1+u_2)m+u_1^2-2u_1u_2+u_2^2+4\}\\
\times \{(u_1u_2+u_1+u_2-1)(u_1u_2-u_1-u_2-1)t^3+6(u_1+u_2)(u_1u_2-1)t^2\\
-3(u_1u_2+u_1+u_2-1)(u_1u_2-u_1-u_2-1)t-2(u_1+u_2)(u_1u_2-1)m^2\\
+(4(u_1+u_2)(u_1u_2-1)t^3-6(u_1u_2+u_1+u_2-1)(u_1u_2-u_1-u_2-1)t^2\\
-12(u_1+u_2)(u_1u_2-1)t+2(u_1u_2+u_1+u_2-1)(u_1u_2-u_1-u_2-1)m\\
-(u_1u_2+u_1+u_2-1)(u_1u_2-u_1-u_2-1)t^3-6(u_1+u_2)(u_1u_2-1)t^2\\
+3(u_1u_2+u_1+u_2-1)(u_1u_2-u_1-u_2-1)t+2(u_1+u_2)(u_1u_2-1)\}\\
\times \{(2(u_1+u_2)(u_1u_2-1)t^3-3(u_1u_2+u_1+u_2-1)(u_1u_2-u_1-u_2-1)t^2\\
-6(u_1+u_2)(u_1u_2-1)t+(u_1u_2+u_1+u_2-1)(u_1u_2-u_1-u_2-1))m^2\\
+(-2(u_1u_2+u_1+u_2-1)(u_1u_2-u_1-u_2-1)t^3-12(u_1+u_2)(u_1u_2-1)t^2\\
+6(u_1u_2+u_1+u_2-1)(u_1u_2-u_1-u_2-1)t+4(u_1+u_2)(u_1u_2-1))m\\
-2(u_1+u_2)(u_1u_2-1)t^3+3(u_1u_2+u_1+u_2-1)(u_1u_2-u_1-u_2-1)t^2\\
+6(u_1+u_2)(u_1u_2-1)t-(u_1u_2+u_1+u_2-1)(u_1u_2-u_1-u_2-1)\}\\
\times \{(t^2+1)^2(m^2+1)(u_1^2+1)(u_2^2+1)\}^{-1}
\end{multline}

The length  of the central diagonal $CF$ is similarly computed since $CF$ is a common side of the two quadrilaterals $ABCF$ and $CDEF$, and equating the circumradii of these two quadrilaterals to $R$ as before, we get two equations for the unknown length $CF$, each equation has four rational roots and the common root gives the diagonal $CF$, which is as follows:
\begin{multline}
CF=\{(u_1+u_2)t-u_1u_2+1\}\{(1-u_1u_2)t-u_1-u_2\}\\
\times \{(u_1+u_2)^2m^2+4(u_1+u_2)m+u_1^2-2u_1u_2+u_2^2+4\}.
\end{multline}

The six minor diagonals of the hexagon, that is, $AC,\;BD,\;AE,\;DF,\;BF$ and $CE$ are also diagonals of the cyclic quadrilaterals $ABCD,\;ADEF$, $ABCF$ and $CDEF$. Since all the four sides of these quadrilaterals are already known, these diagonals are readily found using the formulae \eqref{Brahmaformulae}.   We also note that each minor diagonal divides the hexagon into two parts, one of which is a triangle whose circumradius is $R$ and hence the formula \eqref{Rtrgabc} may also be used to determine the lengths of these diagonals.

The three diagonals $AC$, $BD$ and $BF$ are given by
\begin{equation}
\begin{aligned}
AC&=u_1(t^2+1)(u_2^2+1)(m^2+1)\\
& \quad \quad \{(u_1+u_2)^2m^2+4(u_1+u_2)m+u_1^2-2u_1u_2+u_2^2+4\},\\
BD&=\{(u_1+u_2)((u_1u_2-1)t^2+(2u_1+2u_2)t-u_1u_2+1)m^2\\
& \quad  +((2u_1^2+4u_1u_2-2)t^2+(-4u_1^2u_2+8u_1+4u_2)t-2u_1^2-4u_1u_2+2)m\\
& \quad +(-u_1^2u_2+u_1u_2^2+3u_1+u_2)t^2+(-2u_1^2-4u_1u_2+2u_2^2+4)t\\
& \quad  +u_1^2u_2-u_1u_2^2-3u_1-u_2\}\{(u_1+u_2)((u_1+u_2)t^2+(-2u_1u_2+2)t\\
& \quad-u_1-u_2)m^2+((-2u_1^2u_2+4u_1+2u_2)t^2+(-4u_1^2-8u_1u_2+4)t\\
& \quad+2u_1^2u_2-4u_1-2u_2)m+(-u_1^2-2u_1u_2+u_2^2+2)t^2+(2u_1^2u_2\\
& \quad -2u_1u_2^2-6u_1-2u_2)t+u_1^2+2u_1u_2-u_2^2-2\}(t^2+1)^{-1},\\
BF&=-\{(u_1+u_2)((u_1u_2^2-u_1-2u_2)t^2+(4u_1u_2+2u_2^2-2)t-u_1u_2^2\\
& \quad +u_1+2u_2)m^2+(4(u_1+u_2)(u_1u_2-1)t^2-4(u_1u_2+u_1+u_2-1)\\
& \quad \times (u_1u_2-u_1-u_2-1)t-4(u_1+u_2)(u_1u_2-1)m+\\
& \quad (-u_1^2u_2^2+u_1u_2^3+u_1^2+5u_1u_2-2)t^2-2(u_1+u_2)(2u_1u_2-u_2^2-3)t\\
& \quad+u_1^2u_2^2-u_1u_2^3-u_1^2-5u_1u_2+2\}\{(u_1+u_2)((2u_1u_2+u_2^2-1)t^2\\
& \quad+(-2u_1u_2^2+2u_1+4u_2)t-2u_1u_2-u_2^2+1)m^2+(-2(u_1u_2+u_1\\
& \quad +u_2-1)(u_1u_2-u_1-u_2-1)t^2-8(u_1+u_2)(u_1u_2-1)t\\
& \quad +2(u_1u_2+u_1+u_2-1)(u_1u_2-u_1-u_2-1))m-(u_1+u_2)\\
& \quad \times (2u_1u_2-u_2^2-3)t^2 +(2u_1^2u_2^2-2u_1u_2^3-2u_1^2-10u_1u_2+4)t\\
& \quad +(u_1+u_2)(2u_1u_2-u_2^2-3)\}\{(u_2^2+1)(t^2+1)\}^{-1}
\end{aligned}
\end{equation}

The lengths of the three minor diagonals $DF,\;AE,\;CE$ are given by formulae obtained by interchanging $u_1$ and $u_2$ in the formulae for the lengths of the diagonals $AC,\;BD$ and $BF$ respectively.

Thus the circumradius,  all the 9 diagonals and the area of the cyclic hexagon, whose sides ares given by \eqref{hex}, are all  rational. On suitable scaling, we get a cyclic hexagon whose sides, circumradius, diagonals and area are all integers,

As in the case of the pentagon, it is possible that for certain values of the parameters, we may not get an actual hexagon since  it may not be possible to construct the two quadrilaterals $ABCD$ and $ADEF$ on  opposite sides of the common side $AD$. If $AD$ is twice the circumradius, it becomes a diameter of the circumcircle and the quadrilaterals $ABCD$ and $ADEF$ can always be constructed on opposite sides of $AD$. If $AD$ divides the circumcircle into two unequal segments, and  both the  diagonals of one of the two quadrilaterals are shorter than $AD$ while  one of the diagonals of the second quadrilateral  is longer than $AD$, then also the quadrilaterals $ABCD$ and $ADEF$ can be constructed on opposite sides of $AD$ and we get a cyclic hexagon with the desired properties.  

As a numerical example, if we take the parameters $m,\,t,\,u_1$ and $u_2$ as $5,\,1,\,2$ and $-3$ 
respectively, we get  a cyclic hexagon $ABCDEF$ whose sides $AB,\,BC,\,CD,\,DE$ and $EF$ are 2044, 2880, 3315, 3124, 235, 4420 respectively, the circumradius is 5525/2, the central diagonal $AD$ is 5525 and is thus a diameter of the circumcircle, the area of the heaxgon is 17227230 and the remaining 8 diagonals have lengths $1751561/325,\; 5304,\;4420, \;5133$, $26664/5,\; 3315,\; 4557,$ and  $26167/5$.

As a second numerical example, if we take $(m,\,t,\,u_1,\,u_2)=(2,\,2,\,-8,\,3)$,   
we get a hexagon whose six  sides are $12075,\; 2795,\; 55500,\; 24747,\; 13080$ and $16835$, the circumradius is $60125/2$, the area is 672750078 while the 9 diagonals have lengths $48100,\; 12667468/325,\; 30525,\; 14800,\; 54365,\; 28084,\; 36075$, $ 144943/5$ and $ 533817/13$.

In both the numerical examples given above,  on appropriate scaling, we can easily obtain  a cyclic hexagon with integer sides, diagonals and area.

As in the case of the pentagon, we can try to obtain additional parametrisations of a cyclic hexagon by choosing the  sides of the two quadrilaterals $ABCD$ and $ADEF$  in different ways and imposing the two conditions that the quadrilaterals have the same circumradius and a common side. In some cases, we can obtain a parametrisation for a cyclic rational hexagon but the solution is much more complicated as compared to the solution given above.

\begin{center}
\Large
Acknowledgments
\end{center}
 
I wish to  thank the Harish-Chandra Research Institute, Prayagraj for providing me with all necessary facilities that have helped me to pursue my research work in mathematics.

\medskip

\noindent Postal Address: Ajai Choudhry, 
\newline \hspace{1.05 in}
13/4 A Clay Square,
\newline \hspace{1.05 in} Lucknow - 226001, INDIA.
\newline \noindent  E-mail: ajaic203@yahoo.com

\end{document}